\documentclass{article}
\usepackage{graphicx}
\usepackage{amsmath}
\usepackage{amssymb}
\usepackage{multirow}
\usepackage[a4paper, total={6in, 8in}]{geometry}
\usepackage[colorlinks, citecolor=blue, linkcolor=purple, urlcolor=black]{hyperref}
\usepackage{amsthm}

\usepackage[T1]{fontenc}
\usepackage{xcolor}
\usepackage{tcolorbox}

\tcbset{
    theorembox/.style={
        colback=gray!10,
        colframe=gray!40,
        boxrule=0.5pt,
    }
}

\theoremstyle{definition}
\newtheorem{theorem}{Theorem}
\newtheorem{theorem_1}{Theorem}
\newtheorem{theorem_2}{Theorem}
\title{\textbf{A Complete Decomposition of Stochastic Differential Equations}}
\author{Samuel Duffield \\ 
\small{\href{mailto:sam@normalcomputing.com}{\texttt{sam@normalcomputing.com}}}\\
\small{\textit{Normal Computing}}
}

\begin{document}

\maketitle

\begin{abstract}
    We show that any stochastic differential equation with prescribed time-dependent marginal distributions admits a decomposition into three components: a unique scalar field governing marginal evolution, a symmetric positive-semidefinite diffusion matrix field and a skew-symmetric matrix field.
\end{abstract}

\section{Introduction}

Stochastic differential equations (SDEs) play a fundamental role across numerous scientific and 
engineering disciplines, providing a mathematical framework for modelling systems subject to random 
fluctuations. In mathematical finance, SDEs form the cornerstone of option pricing theory through 
the celebrated Black-Scholes model \cite{black1973pricing}, while in statistical physics and molecular dynamics, they 
describe the evolution of systems towards a thermal equilibrium via Langevin dynamics 
\cite{risken1989fokker}. SDEs have also emerged as a powerful tool in machine learning, 
driving advances in generative modelling through diffusion processes 
\cite{ho2020denoising,song2020score,karras2022elucidating} and enabling efficient sampling algorithms via Markov chain Monte Carlo~\cite{horowitz1991generalized,ma2015complete,duffield2024scalable}.

A fundamental distinction arises between the path distribution $p(x_{[0,T]})$, which characterizes 
the probability measure over entire trajectories $x_{[0,T]} = \{x_t : t \in [0,T]\}$ of the 
stochastic process, and the temporal marginal distributions $p(x,t)$, which describe the 
probability density of the process at each fixed time $t$. While the path distribution captures 
temporal correlations and the full joint statistics of the process across time, the marginal 
distributions $p(x,t)$ only specify the instantaneous state distribution at each moment. In many 
applications, the marginals $p(x,t)$ are of primary interest, as they are directly observable, 
computationally more tractable, and sufficient for characterizing equilibrium states
or generative model outputs. An SDE uniquely determines both, but prescribing only the marginals 
$p(x,t)$ leaves substantial freedom in the path structure, as many different SDEs can share the same 
time-dependent marginals while exhibiting different temporal dependencies.

In this paper, we show that any SDE with prescribed temporal marginal distributions
$p(x,t)$ can be decomposed into three components:
\begin{itemize}
    \item a unique scalar field $\phi(x,t) : \mathbb{R}^d \times [0, T] \to \mathbb{R}$ governing marginal evolution,
    \item a symmetric positive-semidefinite matrix field $D(x,t) : \mathbb{R}^d \times [0, T] \to \mathbb{R}^{d\times d}$,
    \item a skew-symmetric matrix field $Q(x,t) : \mathbb{R}^d \times [0, T] \to \mathbb{R}^{d\times d}$,
\end{itemize}
and that any SDE with prescribed temporal marginal distributions must conform to this
decomposition. Theorem~\ref{theorem:decomp} states the decomposition formally and Theorem~\ref{theorem:unique} states existence and uniqueness of the
scalar field $\phi(x,t)$.

\section{A Complete Decomposition of SDEs}

We now present our complete characterization of SDEs with prescribed temporal marginal 
distributions. The decomposition separates the drift into three interpretable components: 
a scalar field $\phi(x,t)$ that governs the evolution of the marginals as well as a
symmetric positive-semidefinite diffusion matrix $D(x,t)$ and a skew-symmetric matrix $Q(x,t)$ that
are probability preserving. The scalar field $\phi(x,t)$ is linked to $p(x,t)$ through a Poisson 
equation, while $D(x,t)$ and $Q(x,t)$ leave the marginals invariant.

\begin{tcolorbox}[theorembox]
\begin{theorem}\label{theorem:decomp}
    An SDE has temporal marginal distributions $p(x,t)$ (with mild assumptions~\eqref{assumptions_p}) if and only if it has the following form
    \begin{align}
        dx &= \phi(x,t) \nabla_x \log p(x,t)dt +\nabla_x \phi(x,t) dt\nonumber\\
        &\quad+ [D(x,t) + Q(x,t)]\nabla_x \log p(x,t) dt + \nabla_x \cdot [D(x,t) + Q(x,t)]dt \label{eq:decomp}
        \\
        &\quad+ \sqrt{2D(x,t)}dw,\nonumber
        \\
        \partial_t p(x,t) &= - \Delta_x [\phi(x,t)p(x,t)], \label{eq:poisson_eq}
    \end{align}
    for positive-semidefinite $D(x,t) = D(x,t)^\top$ and skew-symmetric $Q(x,t) = -Q(x,t)^\top$.
\end{theorem}
\end{tcolorbox}
\begin{proof}[Proof:~\eqref{eq:decomp} $\implies p(x,t)$]
    An expanded proof is provided in Appendix~\ref{sec:proof_decomp}.

    By~\cite{ma2015complete} the $D(x,t)$ and $Q(x,t)$ terms in~\eqref{eq:decomp} have instantaneous stationary distribution $p(x,t)$ and therefore contribute a zero factor to the Fokker-Planck equation.
    
    This leaves us with
    \begin{equation*}
        dx = \phi(x,t) \nabla_x \log p(x,t) dt + \nabla_x \phi(x,t) dt,
    \end{equation*}
    whose Fokker-Planck equation reads
    \begin{align*}
        \partial_t p(x,t) &= -\nabla_x \cdot [[\phi(x,t) \nabla_x \log p(x,t) + \nabla_x \phi(x,t)]p(x,t)], \\
        &= - \Delta_x [\phi(x,t) p(x,t)],
    \end{align*}
    which is exactly the Poisson equation in \eqref{eq:poisson_eq}.
\end{proof}
\begin{proof}[Proof: \eqref{eq:decomp} $\impliedby p(x,t)$]
    The proof is provided in Appendix~\ref{sec:proof_decomp}. \\ At a high level, it applies
    a Helmholtz decomposition~\cite{glotzl2023helmholtz} to the Fokker-Planck equation
    for a given SDE with prescribed temporal marginal distributions.
\end{proof}

While the diffusion and skew-symmetric components $D(x,t)$ and $Q(x,t)$ can be chosen freely, the 
scalar field $\phi(x,t)$ is uniquely determined by the marginal distributions. This ensures 
that $\phi(x,t)$ carries all the information about how the distribution evolves over time.

\begin{tcolorbox}[theorembox]
\begin{theorem}\label{theorem:unique}
    For given temporal marginal distributions $p(x,t)$ there exists a unique scalar field $\phi(x,t)$ in the decomposition~\eqref{eq:decomp} such that $\lim_{|x| \to \infty} [\phi(x,t)p(x,t)] = 0$.
\end{theorem}
\end{tcolorbox}

\begin{proof}[Proof: Existence] Two constructions for $\phi$ using (i) the fundamental solution to the Laplace operator and (ii) Fourier analysis are provided in Appendix~\ref{sec:proof_exist}. In both cases we restrict to $d \geq 3$ for simplicity.
\end{proof}

\begin{proof}[Proof: Uniqueness]
Suppose $\phi_1(\cdot,t)$ and $\phi_2(\cdot,t)$ both satisfy~\eqref{eq:poisson_eq}, with
$u_i(x,t) := \phi_i(x,t)p(x,t)$. Then
\[
\Delta_x u_1(x,t) = -\partial_t p(x,t) = \Delta_x u_2(x,t),
\]
so $w(x,t):=u_1(x,t)-u_2(x,t)$ satisfies
$
\Delta_x w(x,t)=0.
$
Thus $w(\cdot,t)$ is harmonic on $\mathbb{R}^d$. We assume each $u_i(\cdot,t)$ vanish at infinity and therefore
\[
\lim_{|x|\to\infty} w(x,t)=0.
\]
By the Liouville theorem for harmonic functions~\cite{nelson1961proof}, any harmonic function on $\mathbb{R}^d$ that is bounded—and in particular one that is continuous and vanishes at infinity—must be constant, hence identically zero $w(x,t) = 0 \: \forall x,t$. Therefore $u_1(x,t)=u_2(x,t)$ for all $x$.

Since we assume $p(x,t)>0$ everywhere, it follows that $\phi_1(x,t)=\phi_2(x,t)$ for all $x,t$.
\end{proof}

\section{Related Work}

The most closely related work is the complete recipe for autonomous SDEs and stationary distributions~\cite{ma2015complete}.
This result states that an ergodic autonomous SDE (i.e. an SDE with coefficients independent of time) has a prescribed stationary distribution $\pi(x)$
if and only if it has the following form
\begin{align} \label{eq:complete_recipe}
    dx &= [D(x) + Q(x)]\nabla_x \log \pi(x) dt + \nabla_x \cdot [D(x) + Q(x)]dt
    + \sqrt{2D(x)}dw,
\end{align}
As we will see in Section~\ref{sec:autonomous_sde}, this result can be viewed as a special case of our decomposition~\eqref{eq:decomp}.
The completeness of this result is shown using a Fourier analysis of the Fokker-Planck equation.
A simpler approach using a Helmholtz decomposition (which is also used in our proof of Theorem~\ref{theorem:decomp}) is provided in~\cite{da2023entropy}
in the context of entropy production.
A related decomposition of ergodic autonomous SDEs into symmetric and skew-symmetric components is also provided in \cite{ao2004potential}.

In the context of time-reversed generative differential equations~\cite{anderson1982reverse, haussmann1986time, song2020score}, it is well known that multiple
SDEs (and ODEs) can match the same temporal marginal distributions. Most notably in~\cite{karras2022elucidating} it was shown that Langevin diffusions with different
noise levels can be added without modifying the temporal marginal distributions. Our
decomposition provides a complete characterization of all such SDEs. This is expanded on
in Section~\ref{sec:sde_matching} and Section~\ref{sec:time_reversal}.


The problem of constructing stochastic processes with prescribed marginal distributions is closely related to the Schrödinger bridge problem and, more broadly, to optimal transport.
Given initial and terminal marginals, the Schrödinger bridge selects a diffusion process that interpolates between them while minimizing relative entropy with respect to a reference Brownian motion~\cite{jamison1974reciprocal,follmer1988random,leonard2014survey}. Similarly, there has been work on finding vector fields that solve Liouville's equation \cite{heng2016use} for a given transport path defined by $p(x,t)$. In contrast to this variational perspective, which singles out one distinguished process, our result provides a complete characterization of \emph{all} SDEs consistent with a given family of time-dependent marginals.


\section{Some Corollaries}

\subsection{Autonomous SDEs and Stationary Distributions}\label{sec:autonomous_sde}

We now show that our Theorem~\ref{theorem:decomp} generalizes the autonomous complete recipe from \cite{ma2015complete}.

For an ergodic autonomous SDE with stationary distribution $\pi(x)$ we have
\begin{align*}
    \partial_t p(x,t) = 0.
\end{align*}
This gives Poisson equation \eqref{eq:poisson_eq}
\begin{align*}
    \Delta_x [\phi(x,t)p(x,t)] = 0,
\end{align*}
to which $\phi(x,t) = 0$ is a trivial bounded solution and by Theorem~\ref{theorem:unique} is also unique.

With $\phi(x,t)=0$ and $p(x,t) = \pi(x)$ then \eqref{eq:decomp} exactly matches the form in \eqref{eq:complete_recipe} from \cite{ma2015complete}.

\subsection{SDE Matching}\label{sec:sde_matching}

Suppose we are given an SDE
\begin{align}\label{eq:sde_arb}
    dx = f(x,t)dt + \sqrt{2\Sigma(x,t)}dw,
\end{align}
but not necessarily knowledge of the marginals $p(x,t)$.

By Theorem~\ref{theorem:decomp} the following SDE has the same marginals $p(x,t)$ for any choice of skew-symmetric $Q(x,t)$ and positive-semidefinite $D(x,t)$

\begin{align}
\begin{split}\label{eq:sde_matching}
    dx =& f(x, t)dt + [D(x, t) - \Sigma(x, t) + Q(x, t)]\nabla_x \log p(x, t) dt
\\&+ \nabla_x \cdot [D(x, t) - \Sigma(x, t) + Q(x,t)]dt
+ \sqrt{2D(x,t)} dw,
\end{split}
\end{align}
We can see that~\eqref{eq:sde_matching} matches~\eqref{eq:sde_arb} exactly with $D(x,t) = \Sigma(x,t)$ and $Q(x,t)=0$ and then by Theorem~\ref{theorem:decomp} all SDEs of the form \eqref{eq:sde_matching} have the same marginals for all skew-symmetric $Q(x,t)$ and positive-semidefinite $D(x,t)$.

Note that \eqref{eq:sde_matching} without $Q$ and $D$ terms
\begin{align}\label{eq:sde_match_ode}
    dx = f(x, t)dt - \Sigma(x, t) \nabla_x \log p(x, t) dt
+ \nabla_x \cdot \Sigma(x, t)dt,
\end{align}
has the same marginals but does not necessarily match~\eqref{eq:decomp} without $Q$ and $D$ terms (i.e. the ODE induced by the unique $\phi$) since an alternate $Q$ matrix may already be embedded in \eqref{eq:sde_match_ode}.

\subsection{Time Reversal}\label{sec:time_reversal}

It is well established~\cite{anderson1982reverse,haussmann1986time,cattiaux2023time,kim2025reverse_sde_blog} that an
SDE~\eqref{eq:sde_arb} can be time-reversed to give the following SDE in $y_s = x_{T-s}$
\begin{align}\label{eq:strict_time_reversal}
    dy = - \bar{f}(y, s)ds + 2 \nabla_y \cdot \bar{\Sigma}(y, s)ds + 2 \bar{\Sigma}(y, s) \nabla_y \log \bar{p}(y, s) ds + \sqrt{2\bar{\Sigma}(y, s)}d\bar{w}(s),
\end{align}
where $\bar{f}(y, s) = f(y, T-s)$, $\bar{\Sigma}(y, s) = \Sigma(y, T-s)$, $\bar{p}(y, s) = p(y, T-s)$
and $d\bar{w}$ is a time-reversed Brownian motion.

The SDE~\eqref{eq:strict_time_reversal} is a strict time-reversal in the sense that it
matches the path distribution $p(\bar{y}_{[0, T]}) = p(x_{[0, T]})$ of the original SDE~\eqref{eq:sde_arb}
(here $\bar{y}_s = y_{T-s}$).

In practice~\cite{song2020score}, we are often happy with a weaker time-reversal which only
matches the temporal marginal distributions $p(y, s) = p(x, T-s)$ of the original SDE~\eqref{eq:sde_arb}.
Our Theorem~\ref{theorem:decomp} and Corollary~\ref{sec:sde_matching} applied to~\eqref{eq:strict_time_reversal} provide a
complete characterization of all such time-reversals
\begin{align*}
    dy =& - \bar{f}(y, s)ds
         + [\bar{D}(y, s) + \bar{\Sigma}(y, s) + \bar{Q}(y, s)]\nabla_y \log \bar{p}(y, s) ds
    \\&+ \nabla_y \cdot [\bar{D}(y, s) + \bar{\Sigma}(y, s) + \bar{Q}(y, s)]ds
    + \sqrt{2\bar{D}(y, s)} d\bar{w}(s),
\end{align*}
for any choice of skew-symmetric $\bar{Q}(y, s)$ and positive-semidefinite $\bar{D}(y, s)$.

\subsection{Generative Denoising Models}\label{sec:generative_denoising_models}

A dominant strategy in generative models is to train a parameterised ODE or SDE as the
time-reversal of a specified linear noising SDE with $p(x, 0)$ an empirical
data distribution~\cite{song2020score,huang2021variational,karras2022elucidating}.

Following~\cite{karras2022elucidating}, this noising process can be written as a linear SDE
\begin{align}\label{eq:noising_sde}
    dx = f(t) x dt + g(t)dw,
\end{align}
or equivalently as mollified Gaussian conditionals
\begin{align}\label{eq:noising_gaussian}
    p(x_t \mid x_0) = \mathcal{N}(x_t \mid s(t) x_0, s(t)^2\sigma(t)^2 \mathbb{I}).
\end{align}
See~\cite{karras2022elucidating} for more details including how to convert between
$(f(t), g(t))$ and $(s(t), \sigma(t))$.

Also provided in~\cite{karras2022elucidating} is a recipe for constructing weak
time-reversals of this noising process
\begin{align*}
    dy = - [\partial_s\bar{\sigma}(s)] \bar{\sigma}(s) \nabla_{y} \log \bar{p}(y, s) ds
    + \bar{\beta}(s)\bar{\sigma}(s)^2 \nabla_{y} \log \bar{p}(y, s) ds
     + \sqrt{2 \bar{\beta}(s)} \bar{\sigma}(s) d\bar{w},
\end{align*}
for $\bar{\sigma}(s) = \sigma(T - s)$, $\bar{p}(y, s) = p(y, T-s)$ and any choice of $\bar{\beta}(s) \geq 0$.

Our decomposition elucidates that this recipe is incomplete and a complete recipe of
all valid weak time-reversals is given by
\begin{align}
\begin{split}\label{eq:complete_denoise}
    dy =& - [\partial_s\bar{\sigma}(s)] \bar{\sigma}(s) \nabla_{y} \log \bar{p}(y, s) ds
    \\&+ [\bar{D}(y, s) + \bar{Q}(y, s)]\nabla_y \log \bar{p}(y, s) ds
    + \nabla_y \cdot [\bar{D}(y, s) + \bar{Q}(y, s)]ds
    + \sqrt{2\bar{D}(y, s)} d\bar{w}(s),
\end{split}
\end{align}
for any choice of skew-symmetric $\bar{Q}(y, s)$ and positive-semidefinite $\bar{D}(y, s)$.

The form of~\eqref{eq:complete_denoise} also exactly matches the form in Theorem~\ref{theorem:decomp}
and thus we can conclude that
\begin{align*}
    \phi(y, s) = - [\partial_s\bar{\sigma}(s)] \bar{\sigma}(s),
\end{align*}
is the unique scalar field that satisfies the Poisson equation~\eqref{eq:poisson_eq} for
the weak time-reversal of the noising process (\ref{eq:noising_sde}-\ref{eq:noising_gaussian}).

\section{Discussion}

We have shown that prescribing the full family of temporal marginals $p(x,t)$ essentially fixes one and only one component of the dynamics: the scalar field $\phi(x,t)$ defined by the Poisson equation $\partial_t p(x,t)=-\Delta_x[\phi(x,t)p(x,t)]$ (under mild regularity and decay conditions). All remaining freedom in the path law consistent with the same marginals is captured by a symmetric positive-semidefinite diffusion field $D(x,t)$ and a skew-symmetric field $Q(x,t)$, which are probability-preserving in the sense that they do not affect the instantaneous evolution of the marginals. This yields a complete, constructive characterization of \emph{all} SDEs consistent with a given marginal flow.

Given the generality and completeness of the decomposition, we expect it to motivate many organic research directions; we highlight a few representative avenues. First, the decomposition suggests new optimization and control formulations: one can view $D$ and $Q$ as controllable fields that reshape temporal correlations while keeping $p(x,t)$ fixed, opening the door to objectives for their selection such as minimal entropy production \cite{da2023entropy}, accelerated convergence, or variance reduction for sampling and simulation~\cite{chak2023optimal}. Second, in generative modelling, the characterization provides a complete description of the freedom in weak time-reversals and clarifies how additional noise, non-reversible flows, and state-dependent diffusion can be introduced without affecting the learned marginal path. Specifically, \eqref{eq:complete_denoise} with $Q=D=0$ is often referred to as the probability flow ODE \cite{song2020score, karras2022elucidating} but our formulation points to infinitely many marginal-preserving ODEs dictated by the choice of $Q$. Third, the framework may also be relevant for autonomous SDEs initialized out of stationarity: even when an SDE admits a stationary distribution $\pi(x)$, starting from $p(x,0)\neq \pi(x)$ induces a non-trivial transient marginal flow $p(x,t)$, and the associated (unique) $\phi(x,t)$ isolates the part of the probability current responsible for relaxation towards equilibrium. This viewpoint may be useful for understanding burn-in and designing Markov chain Monte Carlo dynamics that preserve $\pi$ while improving finite-time behaviour, by exploiting the $D$ and $Q$ degrees of freedom without changing the invariant target \cite{ma2015complete,duffield2024scalable}. Fourth, the autonomous recipe \cite{ma2015complete} has been used to devise SDEs that act on an extended space but still preserve the stationary distribution on a subspace (as in Hamiltonian Monte Carlo \cite{betancourt2017conceptual} and underdamped Langevin dynamics \cite{horowitz1991generalized,chak2023optimal}) or via interacting particles \cite{leimkuhler2018ensemble,duncan2023geometry}. Our introduced decomposition enables similar approaches in non-autonomous settings.

Overall, the decomposition provides a simple and general organizing principle for SDEs with prescribed marginals: $\phi$ is the unique signature of marginal evolution, while $D$ and $Q$ parameterise all remaining dynamical freedom.

\section*{Acknowledgements}

We thank Maxwell Aifer, Adrien Corenflos, Samuel Power and Max Welling for discussions and helpful feedback on the manuscript.

\bibliographystyle{plain}
\bibliography{refs}

@article{ma2015complete,
  title={A complete recipe for stochastic gradient {MCMC}},
  author={Ma, Yi-An and Chen, Tianqi and Fox, Emily},
  journal={Advances in Neural Information Processing Systems},
  volume={28},
  year={2015}
}

@article{karras2022elucidating,
  title={Elucidating the design space of diffusion-based generative models},
  author={Karras, Tero and Aittala, Miika and Aila, Timo and Laine, Samuli},
  journal={Advances in Neural Information Processing Systems},
  volume={35},
  pages={26565--26577},
  year={2022}
}

@article{da2023entropy,
  title={The entropy production of stationary diffusions},
  author={Da Costa, Lancelot and Pavliotis, Grigorios A},
  journal={Journal of Physics A: Mathematical and Theoretical},
  volume={56},
  number={36},
  pages={365001},
  year={2023},
  publisher={IOP Publishing}
}

@article{glotzl2023helmholtz,
  title={{Helmholtz} decomposition and potential functions for n-dimensional analytic vector fields},
  author={Gl{\"o}tzl, Erhard and Richters, Oliver},
  journal={Journal of Mathematical Analysis and Applications},
  volume={525},
  number={2},
  pages={127138},
  year={2023},
  publisher={Elsevier}
}

@inproceedings{nelson1961proof,
  title={A proof of {Liouville's} theorem},
  author={Nelson, Edward},
  booktitle={Proc. AMS},
  volume={12},
  pages={995},
  year={1961}
}

@article{black1973pricing,
  title={The pricing of options and corporate liabilities},
  author={Black, Fischer and Scholes, Myron},
  journal={Journal of Political Economy},
  volume={81},
  number={3},
  pages={637--654},
  year={1973},
  publisher={The University of Chicago Press}
}

@incollection{risken1989fokker,
  title={{Fokker}-{Planck} equation},
  author={Risken, Hannes},
  booktitle={The {Fokker}-{Planck} Equation: Methods of Solution and Applications},
  pages={63--95},
  year={1989},
  publisher={Springer}
}

@article{ho2020denoising,
  title={Denoising diffusion probabilistic models},
  author={Ho, Jonathan and Jain, Ajay and Abbeel, Pieter},
  journal={Advances in Neural Information Processing Systems},
  volume={33},
  pages={6840--6851},
  year={2020}
}

@article{duffield2024scalable,
  title={Scalable {Bayesian} learning with posteriors},
  author={Duffield, Samuel and Donatella, Kaelan and Chiu, Johnathan and Klett, Phoebe and Simpson, Daniel},
  journal={arXiv preprint arXiv:2406.00104},
  year={2024}
}

@article{song2020score,
  title={Score-based generative modeling through stochastic differential equations},
  author={Song, Yang and Sohl-Dickstein, Jascha and Kingma, Diederik P and Kumar, Abhishek and Ermon, Stefano and Poole, Ben},
  journal={arXiv preprint arXiv:2011.13456},
  year={2020}
}

@article{anderson1982reverse,
  title={Reverse-time diffusion equation models},
  author={Anderson, Brian DO},
  journal={Stochastic Processes and their Applications},
  volume={12},
  number={3},
  pages={313--326},
  year={1982},
  publisher={Elsevier}
}

@article{haussmann1986time,
  title={Time reversal of diffusions},
  author={Haussmann, Ulrich G and Pardoux, Etienne},
  journal={The Annals of Probability},
  pages={1188--1205},
  year={1986},
  publisher={JSTOR}
}

@article{jamison1974reciprocal,
  title={Reciprocal processes},
  author={Jamison, Benton},
  journal={Zeitschrift f{\"u}r Wahrscheinlichkeitstheorie und Verwandte Gebiete},
  volume={30},
  number={1},
  pages={65--86},
  year={1974}
}

@article{follmer1988random,
  title={Random fields and diffusion processes},
  author={F{\"o}llmer, Hans},
  journal={Lecture Notes in Mathematics},
  volume={1362},
  pages={101--203},
  year={1988},
  publisher={Springer}
}

@article{leonard2014survey,
  title={A survey of the {Schr{\"o}dinger} problem and some of its connections with optimal transport},
  author={L{\'e}onard, Christian},
  journal={Discrete \& Continuous Dynamical Systems - A},
  volume={34},
  number={4},
  pages={1533--1574},
  year={2014}
}

@article{ao2004potential,
  title={Potential in stochastic differential equations: novel construction},
  author={Ao, Ping},
  journal={Journal of Physics A: Mathematical and General},
  volume={37},
  number={3},
  pages={L25},
  year={2004},
  publisher={IOP Publishing}
}

@phdthesis{heng2016use,
  title={On the use of transport and optimal control methods for {Monte Carlo} simulation},
  author={Heng, Jeremy},
  year={2016},
  school={University of Oxford}
}

@article{horowitz1991generalized,
  title={A generalized guided {Monte Carlo} algorithm},
  author={Horowitz, Alan M},
  journal={Physics Letters B},
  volume={268},
  number={2},
  pages={247--252},
  year={1991},
  publisher={Elsevier}
}

@inproceedings{cattiaux2023time,
  title={Time reversal of diffusion processes under a finite entropy condition},
  author={Cattiaux, Patrick and Conforti, Giovanni and Gentil, Ivan and L{\'e}onard, Christian},
  booktitle={Annales de l'Institut {Henri Poincar{\'e}} (B) Probabilit{\'e}s et Statistiques},
  volume={59},
  pages={1844--1881},
  year={2023},
  organization={Institut Henri Poincar{\'e}}
}

@article{huang2021variational,
  title={A variational perspective on diffusion-based generative models and score matching},
  author={Huang, Chin-Wei and Lim, Jae Hyun and Courville, Aaron C},
  journal={Advances in Neural Information Processing Systems},
  volume={34},
  pages={22863--22876},
  year={2021}
}

@misc{kim2025reverse_sde_blog,
  author       = {Ji-Ha Kim},
  title        = {Deriving Reverse-Time Stochastic Differential Equations ({SDEs})},
  year         = {2025},
  month        = {Feb},
  day          = {20},
  url          = {https://jiha-kim.github.io/posts/deriving-reverse-time-stochastic-differential-equations-sdes/},
  organization = {Ji-Ha's Blog},
}

@article{pavliotis2014stochastic,
  title={Stochastic Processes and Applications},
  author={Pavliotis, Grigorios A},
  journal={Texts in Applied Mathematics},
  volume={60},
  year={2014},
  publisher={Springer}
}

@book{stein1970singular,
  title={Singular Integrals and Differentiability Properties of Functions},
  author={Stein, Elias M},
  year={1970},
  publisher={Princeton University Press}
}

@book{folland1995introduction,
  title={Introduction to Partial Differential Equations},
  author={Folland, Gerald B},
  volume={102},
  year={1995},
  publisher={Princeton University Press}
}

@article{duncan2023geometry,
  title={On the geometry of {Stein} variational gradient descent},
  author={Duncan, Andrew and N{\"u}sken, Nikolas and Szpruch, Lukasz},
  journal={Journal of Machine Learning Research},
  volume={24},
  number={56},
  pages={1--39},
  year={2023}
}

@article{betancourt2017conceptual,
  title={A conceptual introduction to {Hamiltonian} {Monte Carlo}},
  author={Betancourt, Michael},
  journal={arXiv preprint arXiv:1701.02434},
  year={2017}
}

@article{leimkuhler2018ensemble,
  title={Ensemble preconditioning for {Markov} chain {Monte Carlo} simulation},
  author={Leimkuhler, Benedict and Matthews, Charles and Weare, Jonathan},
  journal={Statistics and Computing},
  volume={28},
  number={2},
  pages={277--290},
  year={2018},
  publisher={Springer}
}

@article{chak2023optimal,
  title={Optimal friction matrix for underdamped {Langevin} sampling},
  author={Chak, Martin and Kantas, Nikolas and Leli{\`e}vre, Tony and Pavliotis, Grigorios A},
  journal={{ESAIM}: Mathematical Modelling and Numerical Analysis},
  volume={57},
  number={6},
  pages={3335--3371},
  year={2023},
  publisher={EDP Sciences}
}

\begin{appendix}
    

\section{Notation}

\begin{itemize}
\item For a vector field $v : \mathbb{R}^d \to \mathbb{R}^d$, we write $\nabla \cdot v(x) := \sum_{i=1}^d \partial_{x_i} v_i(x)$.

\item For a matrix field $A : \mathbb{R}^d \to \mathbb{R}^{d \times d}$, we define the matrix divergence by
$
[\nabla \cdot A(x)]_i := \sum_{j=1}^d \partial_{x_j} A_{ij}(x)$, with $\nabla \cdot A : \mathbb{R}^d \to \mathbb{R}^d.$

\item For a scalar field $g : \mathbb{R}^d \to \mathbb{R}$, the Laplacian is
$\Delta g(x) := \nabla \cdot \nabla g(x) = \sum_{i=1}^d \partial_{x_i}^2 g(x)$.

\item For $k \in \mathbb{N}$, we write $C^k(\mathbb{R}^d)$ for the space of $k$-times continuously differentiable functions $f : \mathbb{R}^d \to \mathbb{R}$.

\item For $1 \le p <\infty$, we write $L^p(\mathbb{R}^d)$ for the space of measurable functions $f : \mathbb{R}^d \to \mathbb{R}$ such that
$
\|f\|_{L^p} :=
( \int_{\mathbb{R}^d} |f(x)|^p \, dx )^{1/p}
$
is finite.

\item A scalar field $f : \mathbb{R}^d \to \mathbb{R}$ is said to be \emph{bounded} if there exists $M > 0$ such that $|f(x)| \le M$ for all $x \in \mathbb{R}^d$.

\item A function $f : \mathbb{R}^d \to \mathbb{R}^{n}$ is said to \emph{decay sufficiently at infinity} if there exists $\varepsilon > 0$ such that
$
\lim_{|x|\to\infty} |x|^{d+\varepsilon} |f(x)| = 0.
$

\item We say an SDE is \textit{autonomous} if its coefficients are time-independent, i.e. it has the form $dx = b(x)dt + \sqrt{2D(x)}dw$. We say an autonomous SDE is \textit{ergodic} if it has a unique stationary distribution.

\item For a complex number $z = a + b i \in \mathbb{C}$ with $a,b \in \mathbb{R}$ we denote the complex conjugate $\bar{z} = a - bi$.

\item We denote the Fourier transform of a (potentially complex) scalar field $f : \mathbb{R}^d \to \mathbb{C}$ as $\textrm{FT}[f](\xi) = \int e^{-ix^\top\xi}f(x)dx$ with $\textrm{FT}[f]: \mathbb{R}^d \to \mathbb{C}$.

\item We denote the inverse Fourier transform of a (potentially complex) scalar field $g : \mathbb{R}^d \to \mathbb{C}$ as $\textrm{IFT}[g](x) = \frac{1}{(2\pi)^d}\int e^{ix^\top\xi}g(\xi)d\xi$ with $\textrm{IFT}[g]: \mathbb{R}^d \to \mathbb{C}$.

\end{itemize}

\section{Assumptions}

Throughout we make the following mild assumptions on $p(x,t)$:
\begin{subequations}\label{assumptions_p}
\begin{align}
    p(x,t) &> 0 \quad \forall x,t, \\
    p(\cdot,t) &\in C^{2}(\mathbb{R}^d) \quad \forall t, \\
    p(x,\cdot) &\in C^{1}(\mathbb{R}) \quad \forall x, \\
    \int p(x,t) dx &= 1 \quad \forall t, \\
    p(\cdot, t), \nabla_x p(\cdot,t)  \text{ and } \partial_t p(\cdot, t)  &\text{ decay sufficiently at infinity }\forall t.
\end{align}
\end{subequations}

We also note that the decay sufficiently at infinity condition implies $p(\cdot, t), \nabla_x p(\cdot,t), \partial_t p(\cdot, t) \in L^1(\mathbb{R}^d)$. This follows since if $\lim_{|x|\to\infty}|x|^{d+\varepsilon}|f(x)|=0$ for some $\varepsilon>0$, then $|f(x)|\le C|x|^{-d-\varepsilon}$ for $|x|\ge R$, hence
\[
\int_{\mathbb R^d}|f(x)|\,dx \le \int_{|x|\le R}|f(x)|\,dx + C\int_{|x|\ge R}|x|^{-d-\varepsilon}\,dx < \infty,
\]
so $f\in L^1(\mathbb R^d)$.

These assumptions are commonly satisfied in practice and sufficient to justify to ensure validity of the Fokker-Planck equation and suitable integration by parts, see Section 4.1 in \cite{pavliotis2014stochastic}.

For example these assumptions are satisfied for distributions of the form $p(x, t) \propto \exp(-U(x,t))$ where $U(x,t)$ is super-logarithmic in the tails.

\section{Proof of Theorem~\ref{theorem:decomp}}\label{sec:proof_decomp}

We restate Theorem~\ref{theorem:decomp} for the reader's ease.

\begin{tcolorbox}[theorembox]
\begin{theorem_1}
    An SDE has temporal marginal distributions $p(x,t)$ (with mild assumptions~\eqref{assumptions_p}) if and only if it has the following form
    \begin{align}
        dx &= \phi(x,t) \nabla_x \log p(x,t)dt +\nabla_x \phi(x,t) dt\nonumber\\
        &\quad+ [D(x,t) + Q(x,t)]\nabla_x \log p(x,t) dt + \nabla_x \cdot [D(x,t) + Q(x,t)]dt \label{eq:decomp_2}
        \\
        &\quad+ \sqrt{2D(x,t)}dw,\nonumber
        \\
        \partial_t p(x,t) &= - \Delta_x [\phi(x,t)p(x,t)], \label{eq:poisson_eq_2}
    \end{align}
    for symmetric positive-semidefinite $D(x,t) = D(x,t)^\top$ and skew-symmetric $Q(x,t) = -Q(x,t)^\top$.
\end{theorem_1}
\end{tcolorbox}
\begin{proof}[Proof: \eqref{eq:decomp} $\implies p(x,t)$]
    We need to show that the Fokker-Planck equation applied to~\eqref{eq:decomp_2} reduces to the Poisson equation~\eqref{eq:poisson_eq_2}.

    The full Fokker-Planck equation (see also \eqref{eq:J}) reads
    \begin{align*}
        \partial_t p(x,t) &= - \nabla_x \cdot \bigl[[\phi(x,t) \nabla_x \log p(x,t) + \nabla_x \phi(x,t)]p(x,t) \\
        &\quad + [[D(x,t) + Q(x,t)] \nabla_x \log p(x,t)]p(x,t)
        + [\nabla_x \cdot [D(x,t) + Q(x,t)]]p(x,t)\\
        &\quad - \nabla_x \cdot[D(x,t)p(x,t)]\bigr].
    \end{align*}
    Following the result in \cite{ma2015complete} for autonomous SDEs we can show the terms containing $Q$ and $D$ evaluate to zero.

    By linearity of the divergence operator $\nabla \cdot$ we can consider just the $Q$ terms
    \begin{align*}
        Q \text{ terms} &= - \nabla_x \cdot \bigl[Q(x,t) [\nabla_x \log p(x,t)]p(x,t)
        + [\nabla_x \cdot Q(x,t)]p(x,t)\bigr],
        \\
        &=- \nabla_x \cdot \nabla_x \cdot [Q(x,t)p(x,t)]
        \\
        &=0.
    \end{align*}
    Where the second line uses the (reverse) product rule and the third line the skew-symmetry of $Q$.

    Similarly the $D$ terms
    \begin{align*}
        D \text{ terms} &= - \nabla_x \cdot \bigl[D(x,t) [\nabla_x \log p(x,t)]p(x,t)
        + [\nabla_x \cdot D(x,t)]p(x,t)
        - \nabla_x \cdot[D(x,t)p(x,t)] \bigr],
        \\
        &= - \nabla_x \cdot [\nabla_x \cdot D(x,t)p(x,t) - \nabla_x \cdot[D(x,t)p(x,t)]]
        \\
        &=0.
    \end{align*}
    
    This leaves the Fokker-Planck equation as
    \begin{align*}
        \partial_t p(x,t) &= -\nabla_x \cdot [[\phi(x,t) \nabla_x \log p(x,t) + \nabla_x \phi(x,t)]p(x,t)], \\
        &= - \nabla_x \cdot [\nabla_x [\phi(x,t) p(x,t)]],\\
        &= - \Delta_x [\phi(x,t) p(x,t)],
    \end{align*}
    which is exactly the Poisson equation in~\eqref{eq:poisson_eq_2}.
\end{proof}
\begin{proof}[Proof: \eqref{eq:decomp_2} $\impliedby p(x,t)$]
    Now assume we are provided an SDE
    \begin{equation}\label{eq:sde}
        dx = b(x,t)dt + \sqrt{2D(x,t)}dw,
    \end{equation}
    which has temporal marginal distributions $p(x,t)$.
    We now are left to show that there exists a scalar field $\phi(x,t)$ and skew-symmetric matrix field $Q(x,t)$ such that
    \begin{align*}
        b(x,t) &= \phi(x,t) \nabla_x \log p(x,t) +\nabla_x \phi(x,t)
        + [D(x,t) + Q(x,t)]\nabla_x \log p(x,t) + \nabla_x \cdot [D(x,t) + Q(x,t)].
    \end{align*}
    The Fokker-Planck equation for~\eqref{eq:sde} takes the form
    \begin{align}
        \partial_t p(x,t) &= - \nabla_x \cdot J(x,t),  \nonumber\\
        J(x,t) &= b(x,t)p(x,t) - \nabla_x \cdot [D(x,t) p(x,t)]. \label{eq:J}
    \end{align}
    Given assumptions~\ref{assumptions_p} on $p$ and bounded $b, D$ the vector field $J(x,t)$ decays sufficiently at infinity and therefore the Helmholtz decomposition~\cite{glotzl2023helmholtz} states that it decomposes into conservative and divergence-free components
    \begin{align}
        J(x,t) &= \nabla_x[\phi(x,t)p(x,t)] + c(x,t)p(x,t), \label{eq:J_helmholtz}\\
        \nabla_x \cdot [c(x,t)p(x,t)] &= 0. \nonumber
    \end{align}
    By 3.10 in~\cite{glotzl2023helmholtz} we have the vector field $c(x,t)$ satisfies
    \begin{align*}
        \nabla_x \cdot [c(x,t)p(x,t)] = 0 \iff c(x,t) = Q(x,t) \nabla_x \log p(x,t) + \nabla \cdot Q(x,t),
    \end{align*}
    for a skew-symmetric $Q(x,t) = -Q(x,t)^\top$ and $c(x,t)p(x,t) = \nabla_x \cdot [Q(x,t) p(x,t)]$.
    
    Combining the above we get
    \begin{align*}
        b(x,t) p(x,t) &= \nabla_x[\phi(x,t)p(x,t)] + \nabla_x \cdot [Q(x,t) p(x,t)] + \nabla_x \cdot [D(x,t) p(x,t)],
        \\
        \implies b(x,t) &= \phi(x,t) \nabla_x \log p(x,t) +\nabla_x \phi(x,t)\\
        &\quad+ [D(x,t) + Q(x,t)]\nabla_x \log p(x,t) + \nabla_x \cdot [D(x,t) + Q(x,t)].
    \end{align*}
    Assuming $p(x,t) > 0$ everywhere. This now matches~\eqref{eq:decomp_2}.
    
    Finally, the Poisson equation $\eqref{eq:poisson_eq_2}$ is direct from~\eqref{eq:J_helmholtz}
    \begin{align*}
        \partial_t p(x,t) &= -\nabla_x \cdot \nabla_x [\phi(x,t) p(x,t)] - \nabla_x \cdot [c(x,t) p(x,t)],
        \\
        &= - \Delta_x [\phi(x,t) p(x,t)].
    \end{align*}
\end{proof}

\section{Proof of Theorem~\ref{theorem:unique} (Existence)}\label{sec:proof_exist}

We restate Theorem~\ref{theorem:unique} for the reader's ease.

\begin{tcolorbox}[theorembox]
\begin{theorem_2}\label{theorem:unique_2}
    For given temporal marginal distributions $p(x,t)$ there exists a unique scalar field $\phi(x,t)$ in the decomposition~\eqref{eq:decomp} such that $\lim_{|x| \to \infty} [\phi(x,t)p(x,t)] = 0$.
\end{theorem_2}
\end{tcolorbox}

We provide two proofs of the existence of the scalar field $\phi$ solving the Poisson equation~\eqref{eq:poisson_eq} for given distributions $p(x,t)$ satisfying assumptions~\eqref{assumptions_p}. In both cases we assume $d \geq 3$, for a thorough treatment of solutions to Poisson equations including $d=1,2$ see \cite{folland1995introduction}. We also note that by the uniqueness of $\phi$ proved in the main paper the resulting $\phi$s from the two constructions must be equivalent.

\subsection{$\phi(x,t)$ via the Fundamental Solution of the Laplacian}

\begin{proof}[Proof: Existence]
We seek a solution $\phi$ to the Poisson equation
\begin{align*}
    \partial_t p(x,t) &= -\Delta_x u(x,t),
    \\
    u(x,t) &= \phi(x,t)p(x,t).
\end{align*}
Define the fundamental solution of the Laplacian in $\mathbb{R}^d$ (see e.g. Chapter 2 in~\cite{folland1995introduction}) as:
\[
\Phi(x) := 
\frac{|x|^{2-d}}{d(d-2)\omega_d}, \qquad \qquad d\geq3,
\]
where $\omega_d$ denotes the volume of the unit ball in $\mathbb{R}^d$.
For integrable $f \in C^2(\mathbb{R}^d)$ we have that the function
$
g(x) = \int_{\mathbb{R}^d} \Phi(x-y) f(y) dy
$
is well-defined and locally integrable and
\[
\Delta_x \int_{\mathbb{R}^d} \Phi(x-y) f(y) \, dy = -f(x).
\]
So we can set
\begin{align*}
    u(x,t) =  \int_{\mathbb{R}^d} \Phi(x-y) \partial_t p(y,t) dy.
\end{align*}
Noting that our assumptions~\eqref{assumptions_p} imply integrability of $\partial_t p(x,t)$ and $\int \partial_t p(x,t) dx = 0$. Further that $\lim_{|x|\to\infty} u(x,t) = 0$ is also implied through the definition of $\Phi$ and $\partial_t p$ decaying sufficiently at infinity.

This gives
$
\Delta_x u(x,t) = -\partial_t p(x,t).
$
Thus we can set $\phi(x,t) = u(x,t) / p(x,t)$ as $p(x,t) > 0$ everywhere.

\end{proof}

\subsection{$\phi(x,t)$ via Fourier Transform}

\begin{proof}[Proof: Existence]
    
We seek a solution $\phi$ to the Poisson equation
\begin{align*}
    \partial_t p(x,t) &= -\Delta_x u(x,t),
    \\
    u(x,t) &= \phi(x,t)p(x,t).
\end{align*}

Since $\partial_t p$ decays sufficiently at infinity we have the relation between Laplacian and Fourier transform (e.g. Chapter V in \cite{stein1970singular} or Chapter 2-C in~\cite{folland1995introduction})
\begin{align*}
    \textrm{FT}(\Delta_x u)(\xi,t) = - |\xi|^2 \textrm{FT}(u)(\xi,t).
\end{align*}
Thus 
\begin{align*}
    \textrm{FT}(u)(\xi,t) &= -|\xi|^{-2}\textrm{FT}(\Delta_x u)(\xi,t),
    \\
    &= |\xi|^{-2}\textrm{FT}(\partial_t p)(\xi,t) := v(\xi, t),
    \\
    u(x, t) &= \textrm{IFT}[v](x,t).
\end{align*}

Where our assumptions~\eqref{assumptions_p} imply $\int \partial_t p(x,t) dx = 0 \; \forall t$ which in turn implies $\textrm{FT}(\partial_t p)(0,t) = 0$ and therefore $\textrm{FT}(u)(0,t) = 0$ is well defined for $\xi = 0$, further $\partial_t p$ is bounded and continuous so $v(\xi,t)$ is locally integrable near $\xi = 0$ for $d \geq 3$.

Since $\partial_t p(\cdot,t)$ is real-valued, we have the Hermitian symmetry $\textrm{FT}(\partial_t p)(-\xi,t)=\overline{\textrm{FT}(\partial_t p)(\xi,t)}$; as $|\xi|^{-2}$ is real and even, it follows that $v(-\xi,t)=\overline{v(\xi,t)}$, and hence $u(x,t)=\textrm{IFT}[v](x,t)$ is real-valued.

\end{proof}



\end{appendix}

\end{document}